\documentclass[reqno,11pt]{amsart}
\usepackage{graphicx,color}
\newtheorem{theorem}{Theorem}[section]

\newtheorem{lemma}[theorem]{Lemma}

\theoremstyle{definition}

\theoremstyle{remark}

\numberwithin{equation}{section}

\renewcommand{\P}{P}
\newcommand{\E}{E}
\DeclareMathOperator{\diag}{diag}
\DeclareMathOperator*{\argmin}{argmin}
\DeclareMathOperator*{\argmax}{argmax}

\begin{document}

\title[On exponential stability of Wonham filter]
{On exponential stability of Wonham filter}
\author{P. Chigansky}
\address{Dept. Electrical Engineering-Systems,
Tel Aviv University, 69978 Tel Aviv, Israel}
\email{pavelm@eng.tau.ac.il}

\author{R. Liptser}
\address{Department of Electrical Engineering-Systems,
Tel Aviv University, 69978 Tel Aviv, Israel}
\email{liptser@eng.tau.ac.il}


\begin{abstract}
We give elementary proof of a stability result from Atar and
Zeitouni, \cite{AZ},  concerning an exponential
asymptotic ($t\to\infty$) for filtering estimates generated by
wrongly initialized Wonham filter. This proof is based on new
exponential bound having independent interest.
\end{abstract}
\maketitle
\section{Introduction}

In \cite{AZ}, Atar and Zeitouni established  an exponential asymptotic (in
$t\to\infty$) rate of convergence for the distance between differently
initialized filter outputs
 when a signal is
Markov process valued in a finite state space for both the discrete and
continuous time cases (see related topic in \cite{DZ}). In this note, we discuss
only continuous time case\footnote{regarding to discrete
time case we recommend also the recent paper of Le Gland and  Mevel \cite{LM}.}
and deal with exponential stability of the Wonham filter. Following \cite{AZ},
we assume that a signal $X=(X_t)_{t\ge 0}$ is homogeneous Markov process, valued
in a finite alphabet $\mathbb{S}=\{a_1,\ldots,a_d\}$, with transition
intensities matrix $\Lambda$ (with entries $\lambda_{ij}$) and initial
distribution vector $\nu$. An observation $Y=(Y_t)_{t\ge 0}$ is
defined as
\begin{eqnarray*}
Y_t = \int_0^t h(X_s) ds + \sigma W_t
\end{eqnarray*}
where $W=(W_t)_{t\ge 0}$ is a Wiener process independent of $X$, $\sigma$ is a
constant and $h$ is $\mathbb{S}\xrightarrow[]{h}\mathbb{R}$ function.
The Wonham filter, \cite{W} (see also Ch. 9 in\cite{LS}), creates  $\pi_t$
the conditional distribution of $X_t$ given $Y_{[0,t]}=\sigma\{Y_s,s\le t\}$, that is
$
\pi_t(1)=\P(X_t = a_1|Y_{[0,t]}),\ldots,\pi_t(d)=\P(X_t = a_d|Y_{[0,t]}):
$
\begin{equation}
\begin{split}
\label{wf}
&
\pi_0 =\nu
\\
&
d\pi_t = \Lambda^*\pi_tdt + \big(\diag(\pi_t)-\pi_t\pi^*_t\big)h\sigma^{-2}
\big(dY_t-h^*\pi_t dt\big),
\end{split}
\end{equation}
where $^*$, $\diag(\cdot)$, and $h$ are the transposition symbol, scalar matrix, and
vector with entries $h(a_1),\ldots,h(a_d)$ respectively.

It is assumed in \cite{AZ} that
\begin{equation}\label{lambdaij}
\lambda_{ij}>0, \ \forall i\ne j
\end{equation}
and therefore the transition probabilities of Markov process $X$ converge
exponentially fast to the unique and positive invariant distribution. In \cite{AZ},
it is shown that this property is inherited, in some sense,
by $\pi_t$. To emphasize the dependence of $\pi_t$ on $\nu$, below we write
$\pi^\nu_t$ instead of $\pi_t$. Assume $\beta$ the probability distribution
on $\mathbb{S}$, different from $\nu$, is used in \eqref{wf} as the initial condition
and denote $\pi^{\beta\nu}_t$ the solution of corresponding It\^o equation from
\eqref{wf}. For wrong initial distributions $\beta',\beta''$, denote
$\|\pi^{\beta'\nu}_t-\pi^{\beta''\nu}_t\|$  the total variation norm for difference
of wrong conditional distributions.
The following  result  is known from \cite{AZ} (see Theorem 1.6). For any $\beta'$,
$\beta''$
\begin{eqnarray}
\label{rate}
\varlimsup_{t\to\infty}\frac{1}{t}\log
\|\pi^{\beta'\nu}_t-\pi^{\beta''\nu}_t\| \le -2 \min_{p\ne q}\sqrt{\lambda_{pq}
\lambda_{qp}}, \quad P-\text{a.s.}
\end{eqnarray}
This result is very important in applications since it guarantees
that wrongly initialized Wonham filter generates
a wrong conditional distribution converging to correct one
exponentially fast as $t\to\infty$.

The proof of \eqref{rate} in \cite{AZ} (see pp. 46-48 there), being the pioneer,
uses ``heavy'' tools such as Oseledec's multiplicative ergodic theorem,
Birghoff's contraction coefficient, etc. So, it is surprising that
an elementary proof of \eqref{rate} is possible and, moreover,
an exponential estimate for $\|\pi^{\beta\nu}_t-\pi^{\nu}_t\|$, with the same rate
$-2 \min_{p\ne q}\sqrt{\lambda_{pq}\lambda_{qp}}$, is valid for any $t>0$.
\begin{theorem}\label{theo-yat}
Assume $\lambda_{ij}>0$ for any $i\ne j$. Assume $\nu$ is positive and $\beta\sim\nu$.
Then, $P$-a.s. for any $t\ge 0$
\begin{eqnarray}\label{1.3}
\|\pi^{\beta\nu}_t-\pi^{\nu}_t\| \le
 d^2\max_{j\in\mathbb{S}}\frac{d\beta}{d\nu}(j)\max_{j\in\mathbb{S}}\frac{d\nu}{d\beta}
(j)\exp\left(-2
\min_{p\ne q}\sqrt{\lambda_{pq}\lambda_{qp}}t
\right).
\end{eqnarray}
\end{theorem}

\section{Proof of Theorem \ref{theo-yat}}\label{sec-2}

For a notational convenience, write $X^\alpha_t,Y^\alpha_t$ (instead of
$X_t,Y_t$), if $\alpha$ is the distribution of $X_0$.
Then, under $\nu$ being the distribution of $X_0$, write
$$
\pi^\nu_t(i)=P\big(X^\nu_t=a_i|Y^\nu_{[0,t]}\big), \ \ i=1,\ldots,d.
$$
Due to (\ref{wf}), $\pi^{\beta\nu}_t$ is defined as
\begin{equation}
\begin{split}
\label{wf1}
&
\pi^{\beta\nu}_0 =\beta
\\
&
d\pi^{\beta\nu}_t = \Lambda^*\pi^{\beta\nu}_tdt + \big({\sf diag}(\pi^{\beta\nu}_t)-
\pi^{\beta\nu}_t(\pi^{\beta\nu}_t)^*\big)h\sigma^{-2}
\big(dY^\nu_t-h^*\pi^{\beta\nu}_t dt\big).
\end{split}
\end{equation}
It is obvious that
$
\pi^{\beta\nu}_t(i) = \P\big(X^\beta_t=a_i|Y^\beta_{[0,t]}
\big)_{\Big|_{\mbox{$Y^\beta:=Y^\nu$}}},  i=\ 1,\ldots,d.
$

Set
\begin{equation}\label{2.1}
\rho_{ji}(t)=P(X^\nu_0=a_j|Y^\nu_{[0,t]}, X^\nu_t=a_i), \ i,j=1,\ldots,d.
\end{equation}
\subsection{Auxiliary lemmas}
\begin{lemma}\label{lem-2.1}
Under the assumptions of Theorem \ref{theo-yat}, for any $i=1,...,d$
$$
\big|\pi_t^\nu(i)-\pi_t^{\beta\nu}(i)\big|
\le d\max_{j\in\mathbb{S}}\frac{d\beta}{d\nu}(j)
\max_{j\in\mathbb{S}}\frac{d\nu}{d\beta}(j)
\max_{1\le i,j,k\le d}
\big|\rho_{ji}(t)-\rho_{jk}(t)\big|, \ t\ge 0.
$$
\end{lemma}
\begin{proof}
Denote $Q^\nu$ and $Q^\beta$ distributions of $(X^\nu, Y^\nu)$ and
$(X^\beta, Y^\beta)$ respectively. Both $(X^\nu, Y^\nu)$ and
$(X^\beta, Y^\beta)$ are Markov processes, with paths from the  Skorokhod
space $\mathbb{D}_{[0,\infty)}(\mathbb{R}^2)$. They
have the same transition probabilities, $Y^\nu_0\equiv Y^\beta_0 \ (=0)$
and different initial distributions for $X^\nu_0$ and $X^\beta_0$: $\nu$ and $\beta$
respectively. Since $\nu\sim\beta$, the measures $Q^\nu$ and $Q^\beta$ are equivalent and
$$
\frac{dQ^\beta}{d Q^\nu}(X^\nu,Y^\nu)\equiv\frac{d\beta}{d \nu}(X^\nu_0).
$$
We show first that for any $i=1,\ldots,d$ and $t>0$, $Q^\nu$- and $Q^\beta$-a.s.
\begin{eqnarray}
\label{nuform}
\pi^{\beta\nu}_t(i)=\frac
{\sum_{j=1}^d\Big(\frac{d\beta}{d\nu}(j)P\big(X^\nu_0=a_j,
X_t^\nu=a_i|Y^\nu_{[0,t]}\Big)}{E\Big(\frac{d\beta}{d\nu}(X^\nu_0)|Y^\nu_{[0,t]}\Big)}.
\end{eqnarray}

Let $(\mathcal{D}_t)_{t\ge 0}$ be a right continuous filtration on
$\mathbb{D}_{[0,\infty)}(\mathbb{R}^2)$ completed by sets of $Q^\nu$-measure zero.
For fixed $t>0$, let $\psi(x,y)\equiv\psi(y)$, $(x,y)\in\mathbb{D}_{[0,\infty}
(\mathbb{R}^2)$, be $\mathcal{D}_t$ measurable and bounded function. Then
\begin{multline*}
E \psi(Y^\nu) \pi^{\beta\nu}_t(i) E\Big(\frac{d\beta}{d\nu}(X_0^\nu)|Y^\nu_{[0,t]}\Big)
=E \psi(Y^\nu) \pi^{\beta\nu}_t(i)\frac{d\beta}{d\nu}(X_0^\nu)
\\
=E \psi(Y^\nu) \pi^{\beta\nu}_t(i)\frac{dQ^\beta}{dQ^\nu}(X^\nu,Y^\nu)
=E \psi(Y^\beta) \pi^{\beta}_t(i)=\E \psi(Y^\beta) I(X^\beta_t=a_i)
\\
=E \psi(Y^\nu)I(X^\nu_t=a_i)\frac{dQ^\beta}{dQ^\nu}(X^\nu,Y^\nu)
=E \psi(Y^\nu)I(X^\nu_t=a_i)\frac{d\beta}{d\nu}(X^\nu_0)
\\
=E\psi(Y^\nu)\E\Big(I(X^\nu_t=a_i)\frac{d\beta}{d\nu}(X^\nu_0)\big|Y^\nu_{[0,t]}\Big).
\end{multline*}
So, by an arbitrariness of $\psi(Y^\nu)$, we get
$
\pi^{\beta\nu}_t(i)=\frac{E\big(I(X^\nu_t=a_i)\frac{d\beta}{d\nu}(X^\nu_0)|
Y^\nu_{[0,t]}\big)}{E\big(\frac{d\beta}{d\nu}(X_0^\nu)|Y^\nu_{[0,t]}\big)}
$
and it remains to notice that
$$
E\big(I(X^\nu_t=a_i)\frac{d\beta}{d\nu}(X^\nu_0)|
Y^\nu_{[0,t]}\big)=\sum_{j=1}^d\frac{d\beta}{d\nu}(j)P\big(X^\nu_t=a_i,X^\nu_0=a_j|
Y^\nu_{[0,t]}\big).
$$

Now, taking into the consideration \eqref{nuform}, write
\begin{multline*}
\big|\pi_t^\nu(i)-\pi^{\beta\nu}_t(i)\big|
=\Bigg|\pi_t^\nu(i)-\frac
{\sum_{j=1}^d\Big(\frac{d\beta}{d\nu}(j)P\big(X^\nu_0=a_j,
X_t^\nu=a_i|Y^\nu_{[0,t]}\big)\Big)}{E\Big(\frac{d\beta}{d\nu}(X^\nu_0)|Y^\nu_{[0,t]}
\Big)}\Bigg|
\\
=\frac{\Big|\sum_{j=1}^d\frac{d\beta}{d\nu}(j)\Big(
\pi^\nu_t(i)P\big(X^\nu_0=a_j|Y^\nu_{[0,t]}\big)-
P\big(X^\nu_0=a_j,X^\nu_t=a_i|Y^\nu_{[0,t]}\big)\Big)\Big|}
{E\Big(\frac{d\beta}{d\nu}(X^\nu_0)|Y^\nu_{[0,t]}\Big)}.
\end{multline*}
Notice that by the Jensen inequality
$
1\big/E\big(\frac{d\beta}{d\nu}(X^\nu_0)|Y^\nu_{[0,t]}\big)\le
E\big(\frac{d\nu}{d\beta}(X^\nu_0)|Y^\nu_{[0,t]}\big).
$
Consequently,
\begin{multline}\label{formula}
\big|\pi_t^\nu(i)-\pi^{\beta\nu}_t(i)\big|
\le\max_{j\in\mathbb{S}}\frac{d\beta}{d\nu}(j)\max_{j\in\mathbb{S}}\frac{d\nu}{d\beta}(j)
\\
\times
\Bigg|\sum_{j=1}^d\pi^\nu_t(i)\Big(P\big(X^\nu_0=a_j|Y^\nu_{[0,t]}\big)-
P\big(X^\nu_0=a_j,X^\nu_t=a_i,Y^\nu_{[0,t]}\big)\Big)\Bigg|
\\
\le\max_{j\in\mathbb{S}}\frac{d\beta}{d\nu}(j)\max_{j\in\mathbb{S}}\frac{d\nu}{d\beta}(j)
\\
\times\sum_{j=1}^d\pi^\nu_t(i)\Big|P\big(X^\nu_0=a_j|Y^\nu_{[0,t]}\big)-
P\big(X^\nu_0=a_j|X^\nu_t=a_i,Y^\nu_{[0,t]}\big)\Big|
\\
\le\max_{j\in\mathbb{S}}\frac{d\beta}{d\nu}(j)\max_{j\in\mathbb{S}}\frac{d\nu}{d\beta}(j)
\sum_{j=1}^d\Big|P\big(X^\nu_0=a_j|Y^\nu_{[0,t]}\big)-
P\big(X^\nu_0=a_j|X^\nu_t=a_i,Y^\nu_{[0,t]}\big)\Big|
\\
=\max_{j\in\mathbb{S}}\frac{d\beta}{d\nu}(j)\max_{j\in\mathbb{S}}\frac{d\nu}{d\beta}(j)
\sum_{j=1}^d\Big|P\big(X^\nu_0=a_j|Y^\nu_{[0,t]}\big)-\rho_{ji}(t)\Big|.
\end{multline}
Further, an obvious formula
$
P\big(X^\nu_0=j|Y^\nu_{[0,t]}\big)=
\sum_{k=1}^d\pi^\nu_t(k)\rho_{jk}(t),
$
and (\ref{formula}) provide
\begin{multline*}
\big|\pi_t^\nu(i)-\pi^{\beta\nu}_t(i)\big|\le
\max_{j\in\mathbb{S}}\frac{d\beta}{d\nu}(j)\max_{j\in\mathbb{S}}
\frac{d\nu}{d\beta}(j)
\sum_{j=1}^d\Big|\sum_{k=1}^d\pi^\nu_t(k)\rho_{jk}(t)-\rho_{ji}(t)\Big|
\\
=\max_{j\in\mathbb{S}}\frac{d\beta}{d\nu}(j)\max_{j\in\mathbb{S}}\frac{d\nu}{d\beta}(j)
\sum_{j=1}^d\Big|\sum_{k=1}^d\pi^\nu_t(k)\big(\rho_{jk}(t)-\rho_{ji}(t)\big)\Big|
\end{multline*}
and the result.
\end{proof}

\begin{lemma}\label{lem-smooth}
{\rm (Lemma 9.5, \cite{LS}.)}
Assume $\nu$ is positive. Then, for any $j\in\mathbb{S}$,
\begin{equation}\label{smootheq}
\begin{split}
&\rho_{ji}(0)=
  \begin{cases}
    1, & j=i \\
    0, & j\ne i
  \end{cases}
\\
&\frac{d\rho_{ji}(t)}{dt} = \sum_{r\ne i}\frac{\lambda_{ri}
\pi^\nu_t(r)}{\pi^\nu_t(i)}\big(\rho_{jr}(t)-\rho_{ji}(t)\big),
\quad i=1,\ldots,d.
\end{split}
\end{equation}
\end{lemma}

The following lemma plays a key role in proving (\ref{1.3}).
\begin{lemma}\label{lem-key}
Assume $\nu$ is positive and $\lambda_{ij}>0$ for all $i\ne j$. Then
for any $i,j,k$, $P$-a.s.

\begin{equation}
\label{maineq}
|\rho_{ji}(t)-\rho_{jk}(t)|\le \exp\Big(-2\min_{p\ne q}\sqrt{\lambda_{pq}
\lambda_{qp}} t\Big).
\end{equation}
\end{lemma}
\begin{proof}
Set $\rho_\diamond(t)=\min\limits_{i\in\mathbb{S}}\rho_{ji}(t)$ and
$
\rho^\diamond(t)=\max\limits_{i\in\mathbb{S}}\rho_{ji}(t).
$
Obviously, $\rho_\diamond(0)=0$ and $\rho^\diamond(0)=1$. Since processes
$\rho_{ji}(t)$, $i=1,\ldots,d$ have continuously differentiable
paths, piece-wise constant functions\footnote{if $\max$ and $\min$ are attained on
several numbers, the smallest is taken},
$
i^\diamond(t) = \argmax\limits_{i\in\mathbb{S}}\varrho_{ji}(t)
$
and
$
i_\diamond(t) = \argmin\limits_{i\in\mathbb{S}}\rho_{ji}(t)
$
have a finite number of jumps on any finite time interval. Therefore, paths of
$\rho^\diamond(t)$ and $\rho_\diamond(t)$ are
absolutely continuous functions the derivatives of which are defined,
in accordance to \eqref{smootheq}, as
\begin{equation}\label{2.7}
\begin{split}
&
\frac{d\rho_\diamond(t)}{dt}=\sum_{r\ne i_\diamond(t)}
\frac{\lambda_{ri_\diamond(t) }\pi^\nu_t(r)}
{\pi^\nu_t(i_{\diamond}(t))}\big(\rho_{jr}(t)-{\rho}_\diamond(t)\big)
\\
&
\frac{d\rho^\diamond(t)}{dt}=\sum_{r\ne i^\diamond(t)}
\frac{\lambda_{ri^\diamond(t)}\pi^\nu_t(r)}
{\pi^\nu_t(i^\diamond(t))}\big(\rho_{jr}(t)-{\rho}^\diamond(t)\big).
\end{split}
\end{equation}
Notice that for any $r$ we have $\rho_{jr}(t)-\rho_\diamond(t)\ge 0$, and
$\rho_{jr}(t)-\rho^\diamond(t)\le 0$, so that by (\ref{2.7})
\begin{equation*}
\begin{split}
&
\frac{d\rho_\diamond(t)}{dt}\ge\frac{\lambda_{ri_\diamond(t)}\pi^\nu_t(r)}
{\pi^\nu_t(i_\diamond(t))}\big(\rho_{jr}(t)-{\rho}_\diamond(t)\big), \ r\ne i_\diamond(t)
\\
&
\frac{d\rho^\diamond(t)}{dt}\le\frac{\lambda_{ri^\diamond(t)}\pi^\nu_t(r)}
{\pi^\nu_t(i^\diamond(t))}\big(\rho_{jr}(t)-{\rho}^\diamond(t)\big), \ r\ne i^\diamond(t).
\end{split}
\end{equation*}
Consequently, taking $r=i^\diamond(t)$ in the first case and $r=i_\diamond(t)$ for
the second one, we get under $i^\diamond(t)\ne i_\diamond(t)$.
\begin{equation}\label{2.8}
\begin{split}
&
\frac{d\rho_\diamond(t)}{dt}\ge\frac{\lambda_{i^\diamond(t) i_\diamond(t)}
\pi^\nu_t(i^\diamond(t))}
{\pi^\nu_t(i_\diamond(t))}\big(\rho^\diamond(t)-{\rho}_\diamond(t)\big)
\\
&\frac{d\rho^\diamond(t)}{dt}\le-\frac{\lambda_{i_\diamond(t)
i^\diamond(t)}\pi^\nu_t(i_\diamond(t))}
{\pi^\nu_t(i^\diamond(t))}\big(\rho^\diamond(t)-{\rho}_\diamond(t)\big).
\end{split}
\end{equation}
If $i^\diamond(t)=i_\diamond(t)$, then $\rho^\diamond(t)=\rho_\diamond(t)$ and by
(\ref{2.8}) $\frac{d\rho_\diamond(t)}{dt}=0$. Hence, (\ref{2.8}) is always valid.

Set $\triangle_t=\rho^\diamond(t)-\rho_\diamond(t)$ and notice $\triangle_0=1$.
Owing to
\eqref{2.8}, a differential inequality holds true
\begin{eqnarray}
\label{ewe}
\frac{d\triangle_t}{dt}\le -\Big(\frac{\lambda_{i^\diamond(t)i_\diamond(t)}
\pi^\nu_t(i^\diamond(t))}{\pi^\nu_t(i_\diamond(t))}+
\frac{\lambda_{i_\diamond(t)i^\diamond(t)}
\pi^\nu_t(i_\diamond(t))}{\pi^\nu_t(i^\diamond(t))}\Big)\triangle_t.
\end{eqnarray}
Now, an obvious inequality
\begin{equation}\label{rough}
\frac{\lambda_{i^\diamond(t)i_\diamond(t)}\pi_t(i^\diamond(t))}{\pi_t(i_\diamond(t))}
+
\frac{\lambda_{i_\diamond(t)i^\diamond(t)}\pi_t(i_\diamond(t))}{\pi_t(i^\diamond(t))}
\ge \min\limits_{x\ge 0}\Big(\lambda_{i^\diamond(t)i_\diamond(t)}x+
\lambda_{i_\diamond(t)i^\diamond(t)}\frac{1}{x}\Big)
\end{equation}
and
$$
\min\limits_{x\ge 0}\Big(\lambda_{i^\diamond(t)i_\diamond(t)}x+
\lambda_{i_\diamond(t)i^\diamond(t)}\frac{1}{x}\Big)=
2\sqrt{\lambda_{i^\diamond(t)i_\diamond(t)}\lambda_{i_\diamond(t)i^\diamond(t)}}
\ge 2\min_{p\ne q}\sqrt{\lambda_{pq}\lambda_{qp}}
$$
provide
\[
\frac{d\triangle_t}{dt}\le-2\min_{p\ne q}\sqrt{\lambda_{pq}\lambda_{qp}}\triangle_t
\]
and the result. \qed
\end{proof}

\subsection{The proof of (\ref{1.3})}

Since
$
\|\pi^\nu_t-\pi^{\beta\nu}_t\|=\sum_{i=1}^d\big|\pi_t^\nu(i)-\pi_t^{\beta\nu}(i)\big|,
$
by Lemma \ref{lem-2.1} we have
$$
\|\pi_t^\nu-\pi_t^{\beta\nu}\|
\le d^2\max_{j\in\mathbb{S}}\frac{d\beta}{d\nu}\max_{j\in\mathbb{S}}\frac{d\nu}{d\beta}
\max_{1\le i,j,k\le d}
\big|\rho_{ji}(t)-\rho_{jk}(t)\big|, \quad P-\text{a.s.}
$$
and then by Lemma \ref{lem-key} the result. \qed

\section{The proof of (\ref{rate})}

{\bf 1.} Assume $\nu$ is positive and $\beta'\sim\beta''\sim\nu$.
Let $\beta$ denote any of $\beta'$, $\beta''$. Then, due to Theorem \ref{theo-yat}, we
have
\begin{equation}\label{3.1}
\varlimsup_{t\to\infty}\frac{1}{t}\log
\|\pi^{\nu}_t-\pi^{\beta\nu}_t\| \le -2 \min_{p\ne q}\sqrt{\lambda_{pq}
\lambda_{qp}}, \quad P-\text{a.s.}
\end{equation}
Hence, by the triangular inequality
\begin{eqnarray*}
\|\pi^{\beta'\nu}_t-\pi^{\beta''\nu}_t\|&\le& \|\pi^{\beta'\nu}_t-\pi^{\nu}_t\|+
\|\pi^{\beta''\nu}_t-\pi^{\nu}_t\|
\\
&\le& 2\max\Big(\|\pi^{\beta'\nu}_t-\pi^{\nu}_t\|,
\|\pi^{\beta''\nu}_t-\pi^{\nu}_t\|\Big)
\end{eqnarray*}
(\ref{rate}) holds true.

\medskip
{\bf 2.}  It can be readily checked that for any fixed $\delta>0$ the values
$\pi^\nu_\delta$ and $\pi^{\beta\nu}_\delta$, defined by (\ref{wf}) and
(\ref{wf1}) respectively,  have positive entries $Q^\nu$-a.s.
Now, the original filtering problem, under obvious redefinitions, can be translated
into filtering problem on $[\delta,\infty)$ for which \eqref{rate} follows.



\begin{thebibliography}{100}
\bibitem{AZ} R. Atar, O. Zeitouni, Lyapunov exponents for finite state nonlinear filtering,
{\em SIAM J. Control Optim.}, Vol. 35, No 1., pp. 36-55, 1997

\bibitem{DZ}  Delyon, B., Zeitouni, O. Lyapunov exponents for filtering problems in {\em Applied Stochastic Analysis},
ed. M. Davis, R. Elliot, pp. 511-521, Gordon and Breach, New York, 1991

\bibitem{LM} Le Gland, F.,   Mevel, L., Exponential forgetting and geometric ergodicity in Hidden Markov Models,
{\em Math. Control Signals Systems}, Vol. 13, {\bf 1}, pp. 63-93, 2000


\bibitem{LS} Liptser, R., Shiryaev, A., {\em Statistics of Random Processes:
General Theory}, Springer, 2001

\bibitem{W} Wonham, W.M. Some applications of stochastic
differential equations to optimal nonlinear filtering. {\em SIAM J.
Control Optimization}, {\bf 2}, 347--69, 1965.
\end{thebibliography}
\end{document}